\documentclass[11pt]{article}
\newcommand{\documentdate}{10 VII 2024}
\usepackage{a4wide,latexsym,graphicx,amsmath,varioref,mcode}
\usepackage{xcolor}
\title{S2MPJ and {\sf CUTEst} optimization problems \\for Matlab, Python and Julia}

\author{
S. Gratton
  \thanks{Universit\'{e} de Toulouse, INP, IRIT, Toulouse, France. Email:
     serge.gratton@enseeiht.fr. Work partially supported by 3IA Artificial and
     Natural Intelligence Toulouse Institute (ANITI), French "Investing for the Future
     - PIA3" program under the Grant agreement ANR-19-PI3A-0004"}
and Ph. L. Toint
  \thanks{NAXYS, University of Namur, Namur, Belgium. Email: philippe.toint@unamur.be}
}

\newcommand{\beqn}[1]{\begin{equation}\label{#1}}
\newcommand{\eeqn}{\end{equation}}
\newcommand{\req}[1]{(\ref{#1})}
\newcommand{\ms}{\;\;\;\;}

\newcommand{\numsection}[1]{\section{#1}\setcounter{equation}{0}}
\newcommand{\appnumsection}[1]{\section*{#1}\setcounter{equation}{0}
  \renewcommand{\theequation}{A.\arabic{equation}}
  \newcommand{\thetheorem}{A.\arabic{theorem}}
  \renewcommand{\thetable}{A.\arabic{table}}
  \renewcommand{\thefigure}{A.\arabic{figure}}
  \renewcommand{\thesection}{A} }
\newcommand{\ii}[1]{\{ 1, \ldots, #1 \}}

\newcommand{\calE}{{\cal E}} 
\newcommand{\calG}{{\cal G}} 
\renewcommand{\Re}{\hbox{I\hskip -2pt R}}
\newcommand{\smallRe}{\hbox{\footnotesize I\hskip -2pt R}}
\newcommand{\sfrac}[2]{{\scriptstyle \frac{#1}{#2}}}
\newcommand{\half}{\sfrac{1}{2}}
\newcommand{\comment}[1]{}
\renewcommand{\theequation}{\arabic{section}.\arabic{equation}}
\renewcommand{\thefootnote}{(\arabic{footnote})}

\newcommand{\xlower}{x^{\rm low}}
\newcommand{\xupper}{x^{\rm upp}}
\newcommand{\clower}{c^{\rm low}}
\newcommand{\cupper}{c^{\rm upp}}
\newcommand{\cutest}{{\sf CUTEst}}
\newcommand{\thecode}[1]{\begin{quote}\tt #1\end{quote}}

\date{\documentdate}

\begin{document}

\maketitle

\renewcommand{\thefootnote}{\fnsymbol{footnote}}

\begin{abstract}
A new decoder for the SIF test problems of the \cutest\ collection is
described, which produces problem files allowing the computation of
values and derivatives of the objective function and constraints of
most \cutest\ problems directly within ``native'' Matlab, Python or
Julia, without any additional installation or interfacing with MEX
files or Fortran programs. When used with Matlab, the new problem
files optionally support reduced-precision computations.
\end{abstract}


\numsection{Introduction}\label{intro-s}

Published nearly thirty years ago \cite{BongConnGoulToin95} and
updated \cite{GoulOrbaToin03a,GoulOrbaToin15b} since, the {\sf CUTEst}
testing environment and its associated collection of test problems for
continous optimization have been used extensively\footnote{Google
Scholar reports more than 2300 citations, and growing.} by the
mathematical optimization community for the design, testing and
comparison of unconstrained and constrained optimization software. The
test problem collection has indeed grown over time to cover many other
proposals such as the Argonne test set \cite{MoreGarbHill81}, Toint's
collection \cite{Toin83a}, Buckley's problems \cite{Buck89}, COPS
\cite{BondBortMore99}, the Maros and Meszaros quadratic programming
problems \cite{MaroMesz99}, the Hock and Schitkowski collection
\cite{HockSchi81} and the Luk\v{s}an and Vl\v{c}ek test set
\cite{LuksVlce99}, to cite a few of the earlier ones. While the
updates have provided new tools and new problems, the basic design of
the environment has not fundamentally changed since 1995.  Each of the
test problem is described in a file written in ``Standard Input
Format'' (SIF), an (admittedly somewhat obscure) scheme extending the
previous MPS standard \cite{IBM69} for linear programming. To use this
file, it was (and still is) necessary to ``decode'' it using a
decoder/compiler written in Fortran, which produces a (possibly large)
data file ({\tt OUTSDIF.d}) and a small set of Fortran subroutines
encoding the nonlinear parts of the problem. The user then compiles
these Fortran subroutines, loads them with his/her optimization
package and calls the Fortran tools provided in the
\cutest\ environment to evaluate the problem's objective function and 
constraints, possibly with their first and second derivatives.

While this setting has clearly proved its usefulness, its integration
with the evolving computing environment and emerging programming
languages has not always been easy.  An interface with Matlab was
first produced using the ``MEX files'' mechanism available in that
language, but this proved difficult to maintain for all computer
architectures. Other tools such as the {\tt pycutest} interface
\cite{FowkRobeBurm22} for Python, the {\tt matcutest} \cite{Zhan24}
interface to Matlab and CUTEst.jl \cite{SiquOrba23} interface to Julia
have been proposed to make the test problem collection and associated
tools available in today's computing environment, sometimes at the
cost of extra complications. For instance, using the {\tt pycutest}
interface requires the code to be executed in a well-defined
environment running the Fortran decoder, which is typically made easy
by using a ``docker'' container where the original Fortran tools are
precompiled.  The {\tt matcutest} interface does not currently allow
the user to modify the problem parameters (such as dimension or
discretization mesh). Maybe more significantly, all these interfaces
are interfaces with the Fortran SIF decoder and the Fortran {\sf
  CUTEst} evaluation tools.  Although this is fine for computing
environments where Fortran is available or if working in a docker
container is not too restrictive, one must admit that the use of
Fortran has significantly declined since 1995, and that reliance on a
unique tool written in this language could be problematic, in
particular for the continued use of the test problem collection.

The purpose of this paper is to introduce S2MPJ\footnote{For SIF to
  Matlab, Python and Julia.}, a new tool allowing
the computation of values of the objective function and constraints of
most \cutest\ problems, as well as that of their derivatives, directly
within ``native'' Matlab, Python or Julia environments, without any additional
installation or interfacing with MEX files or Fortran programs, using
at most three lines of code.  As it the case for \cutest, S2MPJ also
allows computing the product of the objective function's or
Lagrangian's Hessian or constraints' Jacobian times a user-supplied
vector. When used in the Matlab environment, an option is provided to
use the language variable-precision tools, allowing the user to select
a smaller number of digits in all evaluations.

Our presentation is organized as follows.  After recalling some
details of the structure of the \cutest\ problems in
Section~\ref{GPS-s}, we outline, in Section~\ref{decoder-s}, the main
features and mechanisms allowing S2MPJ to interpret SIF files directly
and describe its formal specifications. We then focus, in
Section~\ref{probfiles-s}, on the description of the newly-decoded
problem files and on how to use them in practice, successively
considering the Matlab, Python and Julia environments. The
reduced-precision option is described in
Section~\ref{mp-s}. Section~\ref{testing-s} presents our methodology
to ensure coherence of the new Matlab, Python and Julia problem files
and evaluation tools with the Fortran-based version.
Section~\ref{distro-s} describes how S2MPJ, problem files and
evaluation tools can be obtained.  Finally, Section~\ref{concl-s}
briefly discusses advantages and drawbacks of our approach and
provides some perpectives for future developments.  However, the
reader eager to use the tool as quickly as possible may directly read
Sections~\ref{how-matlab}, \ref{how-python} and \ref{how-julia}.

\numsection{The group-partially-separable structure of
the \cutest\ \\ problems}\label{GPS-s}

All problems in the \cutest\ collection are specified as
\textit{group-partially-separable} (GPS) problems \cite{ConnGoulToin92}.
Because this is important for understanding what follows, we
remind the reader of what this means. In a GPS problem of the type
encoded in SIF files, one is interested in minimizing a function 
from $\Re^n$ into $\Re$ of the form
\beqn{gps-function}
\min_{x \in \smallRe^n}  f(x)
= \sum_{i \in \calG_{\rm obj}} \frac{F_i[ a_i(x),\omega_i ]}{\sigma_i}
  + \half x^THx,
\eeqn
that is the sum of a quadratic term $\half x^THx$ and of one or more
\textit{objective-function groups} $F_i[ a_i(x), \omega_i ] /
\sigma_i$ for $i$ belonging to some index set $\calG_{\rm obj}$, with
$F_i[.]$  a (possibly nonlinear) univariate \textit{group function}
with parameter(s) $\omega_i$ and
\beqn{group-inside}
a_i(x) = \sum_{j \in \calE_i} w_{ij} f_j(U_j x^{\rm e}_{ij},\tau_{ij})
         - \sum_{j=1}^n \alpha_{ij}\frac{x_j}{\varsigma_j} - \beta_i.
\eeqn
In \req{group-inside}, $\calE_i$ is the index set of \textit{nonlinear
elements} for group $i$, element $j\in\calE_i$ involving the
\textit{element weight} $w_{i,j}$, the (typically nonlinear)
\textit{element function} $f_j$, the (typically low rank and
rectangular) \textit{range matrix} $U_j$, the subvector $x^{\rm
  e}_{ij}$ of \textit{elemental variables} occuring in this element,
the \textit{linear coefficients} $\alpha_{ij}$, the \textit{elemental
  parameter(s)} $\tau_{ij}$, the \textit{group constant} $\beta_i$,
the nonzero \textit{group scaling} $\sigma_i$ and the nonzero
\textit{variable scalings} $\{\varsigma_j\}_{j=1}^n$.  This objective
function may be subject to (possibly infinite) bound constraints of
the form
\beqn{bounds}
\xlower_j \leq x_j \leq \xupper_j \ms (j\in\ii{n})
\eeqn
and to general constraints, each of them involving a single group and
written as
\beqn{cons-format}
\clower_i \leq c_i(x)  \leq \cupper_i   \ms (i\in\ii{m}),
\eeqn
with possibly infinite $\clower_i$ and $\cupper_i$, and
\beqn{cidef}
c_i(x) = \frac{C_i[ a_i(x),\omega_i ]}{\sigma_i}  \ms (i\in\calG_{\rm cons})
\eeqn
where $\calG_{\rm cons}$ is the index set of the groups associated
with constraints, $C_i[\cdot]$ is the univariate constraint's group
function, $\omega_i$ gives its parameter(s) and $\sigma_i$ is the
constraint's scaling. The term $a_i(x)$ also has the form
\req{group-inside} (for $\calE_i$ now the set of nonlinear elements of
group $i \in \calG_{\rm cons}$). Finally, $H$ in \req{gps-function} is
a symmetric $n\times n$ matrix. We let $c(x) = ( c_1(x), \ldots,
c_m(x) )^T$. We also consider the associated Lagrangian function
\[
L(x,y) = f(x) + y^Tc(x),
\]
where $y\in\Re^m$ is the vector of Lagrange multipliers.  As it may
convenient to restrict one's attention to a subset $I = \{i_1, \ldots,
i_p \}$ of the constraints' indeces $\ii{m}$, we also consider the
'$I$-restricted' variants of $c(x)$ and $L(x,y)$ given by
\[
c_I(x) = ( c_{i_1}(x), \ldots, c_{i_p}(x) )^T
\]
and
\beqn{Lagrangian}
L_I(x,y) = f(x) + y_I^Tc_I(x).
\eeqn
where $y_I = ( y_{i_1}(x), \ldots, y_{i_p}(x) )^T$.  Within SIF files,
group functions ($F_i$ and $C_i$) and element functions ($f_j$) are
classified according to named \textit{element types} and \textit{group
  types}, each specifying whether and which parameters are defined for
the type, and, for element types, the number of elemental
variables. When an element $j$ is affected to a group $i$ (that is
$j\in \calE_i$), the name of its type, its weight $w_{ij}$, its
elemental variables $x^{\rm e}_{ij}$ and (optionnally) its elemental
parameters $\tau_{ij}$ are specified. Similarly, the type's name and
(optionally) the group parameters $\omega_i$ are specified for each
group $i$. In addition, SIF allows the definition of \textit{global
  element parameters} and \textit{global group parameters}, that is
parameters whose value is constant across all elements or groups.

The notion of the GPS structure originates in the design of the {\sf
  LANCELOT} package \cite{ConnGoulToin92} where the low rank of the
range matrices $U_j$ was successfully exploited in the ``partitioned
updating'' strategy for quasi-Newton Hessian approximations
\cite{GrieToin82a}. 

Now, given an optimization problem and a vector of variables $x\in
\Re^n$, we may typically be interested in
\begin{itemize}
\item computing the value of the objective function $f(x)$, possibly
  with its
  gradient $g(x)=\nabla_x^1f(x)$ and even
  its Hessian $H(x) = \nabla_x^2f(x)$;
\item computing the vector of constraints' values $c(x$), possibly
  with its Jacobian\footnote{The rows of the Jacobian are the
  transpose of the constraint's gradients.} matrix $J(x) =
  \nabla_x^1c(x)$ and
  even the collection of the constraints' Hessian matrices $\{H_i(x) =
  \nabla_x^2c_i(x)\}_{i=1}^m$;
  \item given a vector $y$ of multipliers, computing the value of the
    Lagrangian function $L(x,y)$, possibly with its gradient (with
    respect to $x$) $\nabla_x^1L(x,y)$ or even its Hessian
    $\nabla_x^2L(x,y)$;
\item given a vector $v\in\Re^n$, computing the product $H(x)v$,
  $J(x)v$ or, given $y$, $\nabla_x^2L(x,y)v$.
\end{itemize}
We may also be interested in performing the same evaluations using a
user-specified ``$I$-restricted'' set of constraints, where $I
\subseteq \ii{m}$. The \cutest\ environment provides Fortran tools for
these evaluations, and so does also our new S2MPJ framework.

It is important to note that the description
\req{gps-function}-\req{Lagrangian} differs from that used
by the Fortran decoder on three points.
\begin{enumerate}
\item By default, S2MPJ numbers the problem's general constraints (in
  $c(x)$, $J(x)$ or $\{H_i(x)\}_{i=1}^m$) starting with ``less or
  equal'' inequality ({\tt <=}) constraints, followed by equality
  ({\tt ==}) constraints, followed by ``larger or equal'' ({\tt >=})
  inequality constraints.  The order in which constraints appear in
  the SIF file is preserved within these three subsets. By contrast,
  the Fortran decoder does not reorder constraints, thus preserving
  the order of appearance in the SIF file. When using S2PMJ, it is
  possible to adhere strictly to the Fortran decoder format when
  decoding the problem's SIF file by specifying the {\tt
    options.keepcorder} flag appropriately (see
  Section~\ref{dec-options}).
\item The format \req{cons-format} in which the constraints are
  described to the user also differs from that used by the Fortran
  decoder, where constraints are described, depending on their type,
  by
  \beqn{For-cons}
  \begin{array}{cclcl}
            r_i \;\leq & c_i(x) &\leq 0 & \ms &  (\mbox{for {\tt <=} constraints})\\
            0 \;=    & c_i(x) &     &         &  (\mbox{for {\tt ==} constraints})\\
            0 \;\leq & c_i(x) & \leq r_i &    &  (\mbox{for {\tt >=} constraints})\\
  \end{array}
  \eeqn
  for \textit{constraint ranges} $r_i$ (optionally) specified in the
  problem's SIF file (see Section 2.1 of the SIF report).  Note that
  $r_i$ must be nonpositive for {\tt <=} constraints, or nonnegative
  for for {\tt >=} ones. This format is, in the authors' view, less
  intuitive, but may be best suited for use with an optimization code
  where inequality constraints are described by transforming them to
  equality ones by the introduction of slack variables. Maintaining
  the Fortran decode format is however possible by setting the
  {\tt options.keepcformat} flag appropriately when decoding the
  problems' SIF file (see Section~\ref{dec-options}).
\item Finally, the SIF format allows the specification of the variable
  scaling factors $\varsigma_j$ in \req{group-inside}. At this time,
  the Fortran decoder ignores them entirely (albeit there are plans to
  provide a tool to pass their values to the user). In S2MPJ, they are
  explicitly taken in to account in that the coefficients of the
  linear terms in each group are adapted according to
  \req{group-inside}. Again, it is possible to adhere to the Fortran
  decoder's choice to ignore them by suitably setting the flag {\tt
    options.pbxscale} when decoding the problem's SIF file, in which
  case the values of the factors $\varsigma_j$ are passed back to the
  user (see Section~\ref{dec-options}), but \ldots full compatibility
  with the Fortran decoder then requires to ignore their value.
\end{enumerate}

\numsection{The S2MPJ framework and its decoder}\label{decoder-s}

The S2MPJ software framework is somewhat simpler than that described
in Section~\ref{intro-s} for the Fortran one. It also starts by
decoding each problem's SIF file (using the {\tt s2mpj.m} Matlab
script), but now produces a single executable ouput file, in Matlab,
Python or Julia, which can then be called for evaluating quantities of
interest directly from the native Matlab, Python or Julia environment,
that is without further interfacing with Fortran.  In order to avoid
making these files longer than necessary, the problem-independent
parts of the process are encapsulated in S2MPJ-supplied
language-dependent libraries ({\tt s2mpjlib.m}, {\tt s2mpjlib.py} and
{\tt s2mpjlib.jl}) which are called while running the problem file.

We now turn to describing the first step of this process (decoding)
and the {\tt s2mpj} script in more detail.

\subsection{Decoder's overview}

   Let {\tt sifpbname} be the name of a problem's SIF file. Broadly
   speaking, this file consists of a ``data'' section, describing the
   structure of the problem as given by
   \req{gps-function}-\req{Lagrangian}, with the exception of the
   definitions for the element functions $f_j(U_j\cdot)$ and the group
   functions $F_i(\cdot)$ and $C_i(\cdot)$. These are then specified
   in the second part of the SIF file, in sections called ELEMENTS and
   GROUPS.
   
   This {\tt sifpbname} file (which must be in the Matlab path) is read
   one line at a time by S2MPJ, and each line is then translated,
   sometimes is a somewhat indirect way, into the corresponding
   Matlab/Python/Julia commands. The decoder
   \begin{itemize}
   \item first sets up the problem data structure, including the value
     of its constant parameters, in the 'setup' action\footnote{That
     is the 'setup' {\tt case} statement in the Matlab output file
     or {\tt \_\_init\_\_ } method in the Python output file or the outer
     {\tt if-then-else}  statement in the Julia output file.};
   \item then includes other actions defining the involved nonlinear
       functions, using information supplied by the Fortran statements in
       the SIF  ELEMENTS and GROUPS sections;
   \item adds a final call to the S2MPJ-supplied libraries {\tt
     s2mpjlib.m} (for Matlab), {\tt s2mpjlib.py} (for Python) or {\tt
     s2mpjlib.jl} (for Julia) that access the previously defined
     actions and data-structures to perform the required evaluation
     tasks.
   \end{itemize}
   As a result, S2PMJ produces an output file called {\tt probname.m},
   {\tt probname.py} or {\tt probname.jl} in the current directory,
   possibly overwriting an existing file with the same name (see
   Section~\ref{dec-options} below to modify the directory where the
   SIF file is found and that where the output file is written).  This
   file is intended for direct calls from Matlab, Python or Julia
   producing values of the objective function, constraints or
   Lagrangian (possibly with their derivatives), as well as products
   of the objective function's or Lagrangian's Hessian or constraints'
   Jacobian times a user-supplied vector (see below). Importantly, the
   produced Matlab/Python file, hereafter called the 'output file', is
   not extensive, meaning that loops are not unrolled (thereby
   maintaining a reasonably compact description).

   The name {\tt probname} of the output file is, in most cases,
   identical to {\tt sifpbname}, but differs from it when the {\tt
     sifpbname} string starts with a digit, in which case it is
   prefixed by {\tt n}, or contains one or more of the characters '+',
   '-', '*' or '/'; these characters are then replaced by {\tt p},
   {\tt m}, {\tt t}, and {\tt d}, respectively. This
   renaming\footnote{For instance, problems {\tt C-RELOAD} and {\tt
     10FOLDTR} are renamed {\tt CmRELOAD} and {\tt n10FOLDTR}.} is
   necessary to allow the output file to be used as a Matlab/Julia
   function or as a Python class.

   Within the 'setup'/\_\_init\_\_ action, instructions are written in
   the output file to define the various problem parameters, as well
   as vectors of bounds on the variables, constraints and objective,
   and the variables/multipliers starting values. Loops in the SIF
   file are directly transformed into loops in the output
   file. Structured sets of the SIF files (such as variables, groups,
   elements) are characterized by the fact that their components may
   be defined using different names and multi-indexes. Because most
   optimization codes only recognize linear structures (a vector of
   variables, a vector of constraints), S2MPJ transforms entities such
   as variables, groups and elements into linear 'flat' unidimensional
   structures. Since this transformation may depend on problem
   internal parameters such as loop limits, themselves depending on
   problem input arguments only known at runtime, the output file uses
   a function\footnote{\tt s2mpj\_ii.} supplied in the {\tt s2mpjlib}
   libraries to compute the relevant indices at runtime.  S2MPJ also
   assign a Matlab/Python/Julia-comptatible name to each of these
   entities. These names are then used in dictionaries associating
   names and values or names and index.

   S2PMJ uses three different structures to pass information between
   the various components of the framework:
   \begin{itemize}
   \item {\tt pbs} is a structure internal to {\tt s1mpj.m} (and
     hence invisible to the user)  used by the code to pass
     information on entities like loops, element's or group's type(s)
     across its different subfunctions;
   \item {\tt pb} is the structure used to pass information on the
     problem to the user once the setup/\_\_init\_\_ action is
     completed (it is described in more detail below);
   \item {\tt pbm} is the structure used to pass information between the
     different actions of the output file (occuring on successive
     calls) and the evaluation functions/methods in the {\tt s2mpjlib}
       librairies. It is visible to the user but should not be
       interfered with (again, it is described more formally below).
   \end{itemize}
   
\noindent
Formally, the {\tt s2mpj} decoder is a Matlab function whose input
parameters are
\begin{description}
\item[\tt sifpbname: ] a string containing the name of the problem to
  be decoded (S2MPJ then reads the {\tt probname.SIF} file for input);
\item[\tt varargin: ] if present, {\tt vararagin\{1\}} allows the
  specification of decoding options (see
  Section~\ref{dec-options}). The use of the suitable option is
  mandatory for decoding SIF files into a Python or a Julia output
  file. Components of {\tt varargin} beyond the first are ignored.
\end{description}
The function {\tt s2mpj.m} has three output values:
\begin{description}
\item[\tt probname: ] a string containing the name of the Matlab,
  Python or Julia output file (without the {\tt .m}, {\tt .py} or {\tt
    .jl} suffix), possibly after renaming (see above);
\item[\tt exitc:] the number of errors which occured before
  termination (or crash). A zero value thus indicates error-free
  execution.
\item[\tt errors:] a cell/list of length {\tt exitc}, whose entries
  contain a brief description of the error(s) found.
\end{description}
Beyond those arguments, the {\tt s2mpj} function produces an output
file called {\tt probname.m}, {\tt probname.py}  or {\tt probname.jl}
in the current directory. 

\subsection{Decoding options}\label{dec-options}

S2MPJ allows the user to specify a number of options, typically to
control file production during execution, save memory, help debugging
or remain as close as possible to the original Fortran SIF
decoder. These options are accessible by suitably setting the first
(and only) variable input argument {\tt varargin\{1\}}.

\vspace*{2mm}
\noindent
If {\tt varargin\{1\}} is a struct then its fields have the following
meaning. 
\begin{description}
\item[\tt varargin\{1\}.language] is is a string specifying the type
  of output produced by S2MPJ and can take different values:
  \begin{description}
  \item[\tt 'matlab': ] a Matlab problem file named {\tt probname.m}
    is written in the current directory, possibly overwriting an
    existing file with the same name. 
  \item[\tt 'python': ] a Python problem file named {\tt probname.py}
    is written in the current directory, possibly overwriting an
    existing file with the same name. 
  \item[\tt 'julia': ] a Julia problem file named {\tt probname.jl}
    is written in the current directory, possibly overwriting an
    existing file with the same name. 
  \item[\tt 'stdma': ] the content of a \textit{potential} Matlab
    problem file is printed on the standard output (no file is produced);
  \item[\tt 'stdpy': ] the content of a \textit{potential} Python
    problem file is printed on the standard output (no file is produced);
  \item[\tt 'stdjl': ] the content of a \textit{potential} Julia
    problem file is printed on the standard output (no file is produced);
  \end{description}
  (default: {\tt 'matlab'})

\item[\tt varargin\{1\}.showsiflines] is a binary flag which is true
  if the SIF data lines must be printed on the standard output. This
  is intended for debugging and requires the option {\tt
    varargin\{1\}.language} to be set to {\tt stdma}, {\tt stdpy} or
  {\tt stdjl}.\\ (default: 0);
 
\item[\tt  varargin\{1\}.getxnames] is a binary flag which is true if
  the names of the variables must be returned in {\tt pb.xnames} on
  exit of the setup action.\\
               (default: 1)

\item[\tt varargin\{1\}.getcnames] is a binary flag which is true if
  the names of the constraints must be returned in
  {\tt pb.cnames} on exit of the setup action.\\
               (default: 1)

\item[\tt varargin\{1\}.getenames] is a binary flag which is true if
  the names of the nonlinear elements $f_j$ must be returned
  in {\tt pbm.elnames} on exit of the setup action.\\
               (default: 0)

\item[\tt varargin\{1\}.getgnames] is a binary flag which is true if
   the names of the groups elements must be returned in
  {\tt pbm.grnames} on exit of the setup action.\\
               (default: 0)

\item[\tt varargin\{1\}.pbxscale] is a binary flag which is true if the
  variable's scaling $\varsigma_j$ (in \req{group-inside}) must be
  provided in {\tt pb.xscale}, instead of being applied internally.\\
               (default: 0)

\item[\tt varargin\{1\}.keepcorder] is a binary flag which is true
  iff, on exit of the setup action, the constraints must not be
  reordered to appear in the order {\tt <=}, followed by {\tt ==} and
  then by {\tt >=}, but should instead appear in the order in which
  they are defined in the SIF file. \\  (default: 0)

\item[\tt varargin\{1\}.keepcformat] is a binary flag which is true
  if, on exit of the setup action, the constraints must be specified
  using the SIF format (i.e.\ specifying ranges and types, see above)
  instead of specifying lower and upper bounds ({\tt clower} and {\tt cupper}) on
  the constraints values. If set, the fields clower and cupper of {\tt pb}
  are replaced by the field {\tt ctypes} of (a cell whose {\tt i}-th
  component contains a string ({\tt '<='}, {\tt '=='}, or {\tt '>=' })
  defining the type of the {\tt i}-th constraint) and, if defined, a
  field {\tt ranges} giving the ranges $r_i$ (see \req{For-cons}) of
  the constraints. \\ (default: 0)

\item[\tt varargin\{1\}.writealtsets] is a binary flag which is true
  if alternative sets of constants, ranges, bounds, starting points or
  objective function bounds must be written (as comments) in the
  output file.\\ (default: 0)

\item[\tt varargin\{1\}.sifcomments] is a binary flag which is true if
  the comments appearing in the SIF file must be repeated in the
  output file.\\ (default: 0)

\item[\tt varargin\{1\}.addinA] is a binary flag which is true if
  repeated entries in the definition of the entries $\alpha_{ij}$ in
  \req{group-inside} must be summed, rather than overwritten.
  Activating this option allows gracefully coping with the SIF error
  of specifying the same entry more than once (as is done by the
  Fortran SIF decoder), albeit at the price of more verbose output
  files.\\ (default:1)

\item[\tt varargin\{1\}.disperrors] is a binary flag which is true if
  the error messages are to be displayed on the standard output as soon as the
  error is detected.\\ (default: 1 )

\item[\tt varargin\{1\}.sifdir] is a string giving the path to the
  directory where the problem's SIF file must be read. (default: in the
  Matlab path)

\item[\tt varargin\{1\}.outdir] is a string giving the path to the
  directory where the output file must be written.\\ (default: '.' )

\item[\tt varargin{1}.redprec] is a binary flag which is true iff
  variable precision support must be provided in the 'setup' section
  of the Matlab output file.\\ (default: 1)
  
\end{description}
\noindent
Not every of the above fields must be defined, each field being tested
for presence and value individually.  If {\tt varargin\{1\}} is not a
Matlab struct or is not present, all options take their default
values.

\subsection{Limitations}\label{limitations-s}
              
   In its current form, the following features are \textit{not}
   supported by S2PMJ:
   \begin{enumerate}
   \item the call to external FORTRAN subroutines;
   \item the D data lines in the data GROUPS section;
   \item the occurrence of blank characters within indeces (probably
     against the standard anyway); 
   \item the FREE FORM version of the SIF file;
   \item... and probably other unforeseen strange SIF and FORTRAN
     constructs.
   \end{enumerate}
   In addition, S2MPJ \textit{does not} provide a comprehensive
   explanation of possible errors in the SIF file (it is supposed to
   adhere to the standard): it may thus crash on errors, or the
   produced output file may also crash. Should this happen, running
   S2MPJ with {\tt varargin\{1\}.showsiflines} $= 1$, usually helps
   spotting the SIF error (or detecting the S2MPJ bug).

   S2MPJ makes the assumption that RANGES $r_i$ (in \req{For-cons})
   are only meaningful for inequality constraints. The SIF standard
   appears to leave open the possibility of ranges for other types of
   groups, but fails to mention what they can/could be used for. Thus
   S2MPJ ignores the ranges of equality constraints or objective
   groups, should them be supplied in the SIF file.

\section{The problem files and how to use them}\label{probfiles-s}
\subsection{The Matlab problem files}
\subsubsection{Matlab problem files: interface specification}

The Matlab output file {\tt probname.m} is is a Matlab function
designed to provide the following user interface. A call of the form
\thecode{varargout = probname( action, varargin )}
produces result(s) stored in {\tt varargout} according
to the following choices of the 'action' argument and the number of
requested outputs, according to the following rules.
\begin{description}      
\item[]\fbox{ \tt [ pb, pbm ] = probname( 'setup',  problem\_parameters ) } \\*[2ex]
  $\bullet$ The fields of the struct {\tt pb} are defined as follows:
  \begin{description}
  \item[\tt pb.name] is a string containing the problem's name
  \item[\tt pb.sifpbname] if present, is a string containing the
    original name of the SIF file before its modification to ensure
    Matlab/Python/Julia compatibility
  \item[\tt pb.n] is an integer giving the problem's number of
    variables
  \item[\tt pb.nob] is an integer giving the number of objective
    groups 
  \item[\tt pb.nle] if present, is an integer giving the number of
    {\tt <=} constraints 
  \item[\tt pb.neq] if present, is an integer giving the number of
    {\tt ==} constraints 
  \item[\tt pb.nge] if present, is an integer giving the number of
    {\tt >=} constraints
  \item[\tt pb.m] is an integer giving the total number of general
    constraints
  \item[\tt pb.lincons] if present, is a integer vector containing the
    indeces of the linear constraints
  \item[\tt pb.pbclass] is the problem's SIF classification
  \item[\tt pb.x0] is a real vector giving the problem's starting
    point
  \item[\tt pb.xlower] is a real vector giving the problem's lower
    bounds on the variables $\xlower_j$ (see \req{bounds})
  \item[\tt pb.xupper] is a real vector giving the problem's upper
    bounds on the variables $\xupper_j$ (see \req{bounds})
  \item[\tt pb.xtype]  is a string whose $i$-th position describes
    the type of the $i$-th variable:
    \begin{description}
    \item[\tt 'r': ] the variable is real
    \item[\tt 'i': ] the variable is integer
    \item[\tt 'b': ] the variable is binary (zero-one)
    \end{description}
    If {\tt xtype} is not a field of {\tt pb}, all variables are
    assumed to be real.
  \item[\tt pb.xscale] if nonempty, is a real vector containing the
    scaling factors $\varsigma_j$ to be applied by the user on the
    occurrences of the variables in the linear terms of the problem's
    groups (see \req{group-inside}). If empty, the scaling
    factors (if any) are applied to the relevant terms internally to the
    decoder without need for further user action (see the {\tt
      pbxscale} option of S2MPJ above)
  \item[\tt pb.y0] if present, is a real vector containing the
    starting values for the constraint's  multipliers
  \item[\tt pb.clower] if present, is the real vector of lower bounds
    on the constraints values $\clower_i$ (see \req{cons-format})
  \item[\tt pb.cupper] if present, is the real vector of upper bounds
    on the constraints values $\cupper_i$ (see \req{cons-format})
  \item[\tt pb.objlower] if present, is a real lower bound of the
    objective function's value
  \item[\tt pb.objupper] if present, is a real upper bound of the
    objective function's value
  \item[\tt pb.xnames] if present, is a cell of strings of length {\tt
    pb.n}  containing the name of the variables 
  \item[\tt pb.cnames] if present, is a cell of strings of length {\tt
    pb.m} containing the name of the constraints.
  \end{description}
   In the 'setup' call to the problem file, {\tt problem\_parameters}
   is a coma-separated list of parameters identified by a \$-{\tt
     PARAMETER} string in the SIF file.  They are assigned in the
   order in which they appear in {\tt varargin}, which is the same as
   that used in the SIF file as well as in the Matlab, Python or Julia
   problem files.  If {\tt varargin} is too short in that it does not
   provide a value for a parameter of the SIF file, the SIF default
   value is used.

   If constraints are present, they are ordered as follows:\\*[1ex]
   \centerline{\begin{tabular}{cl}
    1,\ldots, {\tt pb.nle} &
       : constraints of {\tt <=} type,\\
    {\tt pb.nle+1},\ldots, {\tt pb.nle+pb.neq} &
       : constraints of {\tt ==} type,\\
    {\tt pb.nle+pb.neq+1},\ldots, {\tt pb.nle+pb.neq+pb.nge} &
       : constraints of {\tt >=} type.
   \end{tabular}}

   \noindent
   This ordering also be superseded by setting the {\tt keepcorder}
   flag in the S2MPJ options (see above). If the problem has no
   constraints or bound constraints only, then {\tt lincons}, {\tt
     y0}, {\tt clower}, {\tt cupper} and {\tt cnames} are not fields
   of the {\tt pb} struct returned on setup. {\tt pb.m} is still
   defined in this case, but its value is 0. The fields {\tt objlower}
   and {\tt objupper} may also be missing if no value is provided in
   the SIF file. When specific S2MPJ options are used (again, see
   above), the fields {\tt xnames} and/or {\tt cnames} may be missing
   from the {\tt pb} struct.

   $\bullet$ The {\tt pbm} struct is produced by S2PMJ \textit{for
     information and debugging only. It must not be modified by the
     user}. Its fields are defined as follows:
   \begin{description}
     \item[\tt pbm.name] is a string containing the name of the
       problem,
     \item[\tt pbm.objgrps] if present, is the list of indeces of the
       objective groups, that is $\calG_{\rm obj}$ in \req{gps-function},
     \item[\tt pbm.congrps] if present, is the list of indeces of the
       contraint groups, that is $\calG_{\rm cons}$ in \req{cidef},
     \item[\tt pbm.A] if present, is the real sparse matrix whose
       lines $i$ contains the  coefficients $\alpha_{ij}$ of the
       linear terms in the groups (see \req{group-inside}),
     \item[\tt pbm.Ashape] if present, is an integer vector of length
       2 containing the number of rows of {\tt pbm.A} and its number of
       columns (for Python only), 
     \item[\tt pbm.gconst] if present, is a real vector containing the
       group's constants $\beta_i$ of \req{group-inside},
     \item[\tt pbm.H] if present, is a symmetric sparse real matix
       containing the Hessian matrix $H$ of \req{gps-function},
     \item[\tt pbm.enames] if present, is a cell of strings containing
       the names of the nonlinear elements,
     \item[\tt pbm.elftype] if present, is a cell of strings
       containing the names of the element's types,
     \item[\tt pbm.elvar] if present, is a cell of integer vectors
       containing the indeces of the elemental variables for each
       element,
     \item[\tt pbm.elpar] if present, is a cell of real vectors
       containing the values of the elemental parameters for each
       element,
     \item[\tt pbm.gscale] if present, is a real vector containing the
       group's scaling factors $\sigma_i$ in \req{gps-function},
     \item[\tt pbm.grnames] if present, is a cell of strings
       containing the group's names,
     \item[\tt pbm.grftype] if present, is a cell of strings
       containing the group's types,
     \item[\tt pbm.grelt] if present, is a cell of integer vectors
       containing the indeces of the nonlinear elements occuring in
       each group,
     \item[\tt pbm.grelw] if present, is a cell of real vectors
       containing the weights of the nonlinear elements occuring in
       each group, that is the $w_{ij}$ in \req{group-inside},
     \item[\tt pbm.grpar] if present, is a cell of real vectors
       containing the values of the parameters $\omega_i$ of the group
       functions,
     \item[\tt pbm.efpar] if present, is a real vector containing the
       values of the global element parameters,
     \item[\tt pbm.gfpar] if present, is a real vector containing the
       values of the global group parameters,
     \item[\tt pbm.ndigs] if present is an integer giving the number of
       digits requested for reduced-precision evaluations (see
       Section~\ref{mp-s}). 
     \end{description}
     Depending on the problem's nature, some fields may be missing
     from the {\tt pbm} struct, but at least one of {\tt objgrps} or
     {\tt congrps} must be present and non-empty. The presence of the
     fields {\tt enames} and {\tt gnames} depends on the {\tt
       getenames} and {\tt getgnames} S2MPJ options (see
     Section~\ref{dec-options}).

     \fbox{\parbox{0.99\linewidth}{A call to
     {\tt probname( 'setup', ...)} \textit{must precede} any call to
     {\tt probname} with other actions (in order to setup the
     problem's data structure).}} 

     \vspace*{2mm}
\item[]\fbox{ \tt fx = probname( 'fx',  x ) } \\ \vspace*{-3mm}
  \begin{description}
    \item[\tt fx] is $f(x)$, the value of the objective
      function at $x$, 
  \end{description}
  
\item[]\fbox{ \tt [ fx, gx ] = probname( 'fgx',  x ) } \\ \vspace*{-3mm}
  \begin{description}
    \item[\tt fx] is $f(x)$, the value of the objective
      function at $x$, 
    \item[\tt gx] is $\nabla_x^1f(x)$, the value of the objective
      function's gradient at $x$, 
  \end{description}
  Note that the call { \tt [ fx, gx ] = probname( 'fx',  x ) }
  produces the same results.

\item[]\fbox{ \tt [ fx, gx, Hx ] = probname( 'fgHx',  x ) } \\ \vspace*{-3mm}
  \begin{description}
    \item[\tt varargout\{1\}] is $f(x)$, the value of the objective
      function at $x$, 
    \item[\tt varargout\{2\}] is $\nabla_x^1f(x)$, the value of the
      objective function's gradient at $x$,
    \item[\tt varargout\{3\}] is $\nabla_x^2f(x)$, the value of the
      objective function's Hessian at $x$,
  \end{description}
  Note that the call { \tt [ fx, gx, Hx ] = probname( 'fx',  x ) }
  produces the same results.

\item[]\fbox{ \tt cx = probname( 'cx',  x ) } \\ \vspace*{-3mm}
  \begin{description}
    \item[\tt cx] is $c(x)$, the vector of constraint's
      values at $x$,  
  \end{description}

\item[]\fbox{ \tt [ cx, Jx ] = probname( 'cJx',  x ) } \\ \vspace*{-3mm}
  \begin{description}
    \item[\tt cx] is $c(x)$, the vector of constraint's
      values at $x$,  
    \item[\tt Jx] is the matrix $J(x)$, the constraint's
      Jacobian at $x$,  
    \end{description}
  Note that the call { \tt [ fx, gx ] = probname( 'cx',  x ) }
  produces the same results.

\item[]\fbox{ \tt [ cx, Jx, Hix ] = probname( 'cJHx',  x ) } \\ \vspace*{-3mm}
  \begin{description}
    \item[\tt cx] is $c(x)$, the vector of constraint's
      values at $x$,  
    \item[\tt Jx] is the matrix $J(x)$, the constraint's
    \item[\tt Hix] is a cell whose $i$-th entry is
      $\nabla_x^2c_i(x)$, the Hessian of the $i$-th constraint at $x$,
  \end{description}
  Note that the call { \tt [ fx, gx, Hix ] = probname( 'cx',  x ) }
  produces the same results.

\item[]\fbox{ \tt cIx = probname( 'cIx',  x, I ) } \\ \vspace*{-3mm}
  \begin{description}
    \item[\tt cIx] is $c_I(x)$, the vector of the
      $I$-restricted constraint's values at $x$,  
  \end{description}

\item[]\fbox{ \tt [ cIx, JIx ] = probname( 'cIJx',  x, I ) } \\ \vspace*{-3mm}
  \begin{description}
    \item[\tt cIx] is $c_I(x)$, the vector of the
      $I$-restricted constraint's values at $x$,  
    \item[\tt JIx] is the matrix $J_I(x)$, the Jacobian of
      the $I$-restricted constraints at $x$,  
  \end{description}

\item[]\fbox{ \tt [ cIx, JIx, HIx ] = probname( 'cIJHx',  x, I ) } \\ \vspace*{-3mm}
  \begin{description}
    \item[\tt cIx] is $c_I(x)$, the vector of the
      $I$-restricted constraint's values at $x$,  
    \item[\tt JIx] is the matrix $J_I(x)$, the Jacobian of
      the $I$-restricted constraints at $x$,  
    \item[\tt HIx] is a cell whose $i$-th entry is
      $\nabla_x^2c_\ell(x)$, the Hessian of the $\ell$-th constraint at $x$
      where $\ell$ is the $i$-th entry in $I$,
  \end{description}

\item[]\fbox{ \tt Hxv = probname( 'fHxv',  x, v ) } \\ \vspace*{-3mm}
  \begin{description}
    \item[\tt Hxv] is the vector $\nabla_x^2f(x)v$, the
      product of the objective function's Hessian at $x$ times the
      vector $v$,
  \end{description}
  
\item[]\fbox{ \tt Jxv = probname( 'cJxv',  x, v ) } \\ \vspace*{-3mm}
  \begin{description}
    \item[\tt Jxv] is the vector $J(x)v$, the
      product of the constraint's Jacobian at $x$ times the
      vector $v$,
  \end{description}
  
\item[]\fbox{ \tt JIxv = probname( 'cIJxv',  x, v ) } \\ \vspace*{-3mm}
  \begin{description}
    \item[\tt JIxv] is the vector $J_I(x)v$, the
      product of the constraint's $I$-restricted Jacobian at $x$ times the
      vector $v$,
  \end{description}
  
\item[]\fbox{ \tt Lxy = probname( 'Lxy',  x, y ) } \\ \vspace*{-3mm}
  \begin{description}
    \item[\tt Lxy] is $L(x,y)$, the value of the problem's Lagrangian
      at $(x,y)$, 
  \end{description}
  
\item[]\fbox{ \tt [ Lxy, gLxy ] = probname( 'Lgxy',  x, y ) }\\ \vspace*{-3mm}
  \begin{description}
    \item[\tt Lxy] is $L(x,y)$, the value of the problem's Lagrangian
      at $(x,y)$, 
    \item[\tt gLxy] is $\nabla_x^1L(x,y)$, the value of the Lagrangian's
      gradient with respect to $x$ taken at $(x,y)$, 
  \end{description}
  Note that the call { \tt [ Lxy, Lgxy ] = probname( 'Lxy',  x, y ) }
  produces the same results.

\item[]\fbox{ \tt [ Lxy, gLxy, HLxy ] = probname( 'LgHxy',  x, y ) }\\ \vspace*{-3mm}
  \begin{description}
    \item[\tt Lxy] is $L(x,y)$, the value of the problem's Lagrangian
      at $(x,y)$, 
    \item[\tt gLxy] is $\nabla_x^1L(x,y)$, the value of the Lagrangian's
      gradient with respect to $x$ taken at $(x,y)$, 
    \item[\tt HLxy] is $\nabla_x^2L(x,y)$, the value of the Lagrangian's
      Hessian with respect to $x$ taken at $(x,y)$, 
  \end{description}
  Note that the call { \tt [ Lxy, gLxy, HLxy ] = probname( 'Lxy', x, y ) }
  produces the same results.

\item[]\fbox{ \tt LIxy = probname( 'LIxy',  x, y, I ) } \\ \vspace*{-3mm}
  \begin{description}
    \item[\tt LIxy] is $L_I(x,y)$, the value of the problem's
      $I$-restricted Lagrangian at $(x,y)$, 
  \end{description}
  
\item[]\fbox{ \tt [ LIxy, gLIxy ] = probname( 'LIgxy',  x, y, I ) }\\ \vspace*{-3mm}
  \begin{description}
    \item[\tt LIxy] is $L_I(x,y)$, the value of the problem's $I$-restricted Lagrangian
      at $(x,y)$, 
    \item[\tt gLIxy] is $\nabla_x^1L_I(x,y)$, the value of the $I$-restricted Lagrangian's
      gradient with respect to $x$ taken at $(x,y)$, 
  \end{description}
  Note that the call { \tt [ LIxy, gLIxy ] = probname( 'Lxy',  x, y, I ) }
  produces the same results.

\item[]\fbox{ \tt [ LIxy, gLIxy, HLIxy ] = probname( 'LIgHxy',  x, y, I ) }\\ \vspace*{-3mm}
  \begin{description}
    \item[\tt LIxy] is $L_I(x,y)$, the value of the problem's
      $I$-restricted Lagrangian at $(x,y)$,
    \item[\tt gLxy] is $\nabla_x^1L_I(x,y)$, the value of the
      $I$-restricted Lagrangian's gradient with respect to $x$ taken
      at $(x,y)$,
    \item[\tt HLxy] is $\nabla_x^2L_I(x,y)$, the value of the
      $I$-restricted Lagrangian's Hessian with respect to $x$ taken at
      $(x,y)$,
  \end{description}
  Note that the call { \tt [ LIxy, gLIxy, HLIxy ] = probname( 'LIxy',  x, y ) }
  produces the same results.

\item[]\fbox{ \tt HLxyv = probname( 'LHxyv',  x, y, v )}\\ \vspace*{-3mm}
    \begin{description}
    \item[\tt HLxyv] is the vector $\nabla_x^2L(x,y)v$, the product
      the problem's Lagrangian Hessian with respect to $x$ taken at
      $(x,y)$ times the vector $v$,
    \end{description}
      
\item[]\fbox{ \tt HLIxyv = probname( 'LIHxyv',  x, y, v, I)}\\ \vspace*{-3mm}
    \begin{description}
    \item[\tt HLIxyv] is the vector $\nabla_x^2L_I(x,y)v$, the product
      the problem's $I$-restricted Lagrangian Hessian with respect to
      $x$ taken at $(x,y)$ times the vector $v$.
    \end{description}
      
\end{description}
The actions {\tt cx}, {\tt cIx}, {\tt Jxv}, {\tt JIxv}, {\tt Lxy},
{\tt LIxy}, {\tt HLxyv} or {\tt HLIxyv} are of course meaningless when the
problem is unconstrained or only has bound constraints.

\subsubsection{Matlab problem files: how to use them}\label{how-matlab}

The use of S2MPJ in Matlab is typically as follows.
\begin{enumerate}
  \item a SIF problem ({\tt PROBNAME.SIF}, say) is decoded by the
    Matlab command 
    \thecode{ PROBLEM = s2mpj( 'SIFPROBLEM' );}
    or
    \thecode{PROBLEM = s2mpj( 'SIFPROBLEM', options );}
    where options is a struct whose content is described above
    in Section \ref{dec-options}. This assumes that the
    {\tt PROBNAME.SIF} file is in the Matlab path. A file {\tt
      PROBLEM.m} is then produced in the current directory
    (potentially overwriting an existing one with the same name).
\item The problem data structure is then setup, along with the
  starting points, bounds and other components of the {\tt pb} and
  {\tt pbm} structs (see above) by issuing the command
  \thecode{pb = PROBLEM( 'setup', args\{:\} )}
   where {\tt args\{:\}} is an optional comma-separated list of
   problem-dependent arguments (such as problem's dimension, for
   instance) identified by the string \$-{\tt PARAMETER} in the SIF
   file. If {\tt args\{:\}} is missing, setup is performed using the
   SIF file defaults.
 \item The value(s) of the problem functions (objective and
  constraints) at a vector $x$, together with values of their 
  first and second derivatives (if requested) are then computed by
  issuing one of the commands
  \thecode{fx = PROBLEM( 'fx', x );}
  \thecode{[ fx, gx ] = PROBLEM( 'fgx', x );}
  \thecode{[ fx, gx, Hx ]  = PROBLEM( 'fgHx', x );}
  for the objective function and (if constraints are present)
  \thecode{cx = PROBLEM( 'cx', x );}
  \thecode{[ cx, Jx ] = PROBLEM( 'cJx', x );}
  \thecode{[ cx, Jx, Hx ]  = PROBLEM( 'cJHx', x ); }
  for the constraints.  In this case, it is possible to restrict one's
  attention to a subset $I$ of the constraints (i.e.\ using $c_I(x)$,
  the '$I$-restricted' version of $c(x)$) by using
  \thecode{cIx = PROBLEM( 'cIx', x, I );}
  \thecode{[ cIx, JIx ] = PROBLEM( 'cIJx', x, I );}
  \thecode{[ cIx, JIx, HIx ] = PROBLEM( 'cIJHx', x, I );}
  The product of the objective function's Hessian (at $x$) times a
  user-supplied vector $v$ can be obtained by the command
  \thecode{Hxv = PROBLEM( 'fHxv', xx, v );}
   while the product of the (potentially $I$-restricted) constraints'
   Jacobian (at $x$) times $v$ is computed by issuing one of the
   commands
   \thecode{Jxv  = PROBLEM( 'cJxv' , x, v );}
   \thecode{JIxv = PROBLEM( 'cIJxv', x, v, I );}
   The value (and derivatives) of the Lagrangian function
   $L(x,y)$ at $(x,y)$  is obtained by one of the commands 
   \thecode{Lxy = PROBLEM( 'Lxy', x, y );}
   \thecode{[ Lxy, Lgxy ] = PROBLEM( 'Lgxy', x, y );}
   \thecode{[ Lxy, Lgxy, LHxy ] = PROBLEM( 'LgHxy', x, y );}
   while the product of the Lagrangian's Hessian times a vector $v$
   can be computed with the command
   \thecode{HLxyv = PROBLEM( 'HLxyv', x, y, v );}
   Finally, and as above for constraints, the Lagrangian may be
   $I$-restricted in the commands 
   \thecode{[LIxy, gLIxy, HLIxy ] = PROBLEM( 'LIxy', x, y, I );}
   \thecode{HLIxyv = PROBLEM( 'LIHxyv' , x, y, v, I );}
   Because {\tt PROBNAME.m} uses a persistent structure to pass problem
   structure across its various actions, such calls are valid as long
   as the Matlab variable space is not cleared and as long as another
   {\tt OTHERPROBLEM.m} file is not called. 
\end{enumerate}
    
\subsection{The Python problem files}
\subsubsection{Python problem files: interface specification}

   The Python output file describes a Python class whose name is that of the
   problem and which is derived from the parent {\tt CUTEst\_problem}
   class described in the {\tt s2pmjlib.py} file. Its inherits the
   evaluation methods of this parent class, which are organized in a
   manner similar to actions for the Matlab output file, except that
   the Matlab 'setup' action is replaced by a simple call to the class
   with the problem parameters, the setting up of the problem
   structure(s) being then performed by the {\tt \_\_init\_\_} method
   of the class. It also inherits all fields described above for
   the {\tt pb} and {\tt pbm} structs in Matlab. So, if the problem is
   given by the {\tt PROBNAME} class, setting it up is achieved by the
   call
   \thecode{problem =  probname(args[:]}
   where {\tt args[:]} is an optional comma-separated list of
   problem-dependent arguments.  The problem's
   starting point is then given by {\tt problem.x0}, while the
   evaluation task corresponding to the Matlab call {\tt [outputs] =
     probname( action, x, other\_args );} is performed by the Python call
   \thecode{ouputs = problem.action( x, other\_args )} for all
   actions (except 'setup') listed above. The real vectors resulting
   from the various evaluations or being present in the fields of {\tt
     problem} are column-oriented numpy.ndarrays.  Real
   matrices are sparse.csr matrices, or lists of sparse.csr matrices
   (for the Hessians of the constraints).

\subsubsection{Python problem files: how to use them}\label{how-python}

   Again the use of the Python output file is very similar to that of
   the Matlab file. The various steps described above for the the use
   of the Matlab output file are now (given the {\tt PROBNAME.SIF} file) as
   follows.  The decoding of the SIF file and production of the
   {\tt PROBNAME.py} file containing the {\tt PROBNAME} class is
   performed by the call
   \thecode{PROBNAME = s2mpj( 'SIFPBNAME', inpy )}
   where {\tt inpy} is a struct describing the S2MPJ options (see
   Section\ref{dec-options}) with {\tt inpy.language = 'python'} (other
   options may be specified if desired.)  After importing the
   functions of the {\tt PROBNAME} module by 
   \thecode{from PROBNAME import *}
   the setup of the problem is now performed by a call of the form
   \thecode{PROBLEM = PROBNAME( arg[:] )}
   where {\tt args[:]} is a comma separated list of problem-dependent
   arguments. The subsequent evaluation tasks are then obtained by
   calling one (or more) of
   \begin{quote}\tt
   fx          = PROBLEM.fx( x )\\
   fx, gx      = PROBLEM.fgx( x )\\
   x, gx, Hx   = PROBLEM.fgHx( x )\\
   Hxv         = PROBLEM.fHxv( x, v )\\
   cx          = PROBLEM.cx( x )\\
   cx, Jx      = PROBLEM.cJx( x )\\
   cx, Jx, HXs = PROBLEM.cJHx( x )\\
   Jxv         = PROBLEM.cJxv( x, v )\\
   cIx         = PROBLEM.cIx( x, I )\\
   cIx, cIJx   = PROBLEM.cIJx( x, I )\\
   cIx, cIJx, cIJHx = PROBLEM.cIJHx( x, I )\\
   cIHv        = PROBLEM.cIHxv( x, v, I )\\
   Lxxy        = PROBLEM.Lxy( x, y )\\
   Lxy, Lgxy   = PROBLEM.Lgxy( x, y )\\
   Lxy, gLxy, HLxy = PROBLEM.LgHxy( x, y )\\
   LHxyv       = PROBLEM.LHxyv( x, y, v )\\
   LIxy        = PROBLEM.LIxy( x, y, I )\\
   LIxy, gLIxy = PROBLEM.LIgxy( x, y, I )\\
   LIxy, gLIxy, HLIxy = PROBLEM.LIgHxy( x, y, I )\\
   LIHxyv      = PROBLEM.LIHxyv( x, y, v, I )
   \end{quote}

\subsection{The Julia problem files}
\subsubsection{Julia problem files: interface specification}
       
   The Julia output file is a cross between the Matlab and Python
   ones. It is a Julia function whose interface is identical to that
   of the Matlab output, \textit{except} that the second argument in
   the call must be the {\tt pbm} struct. This struct now has an
   additional field {\tt call} so that {\tt pbm.call} contains the
   Julia call to the output file after execution with the 'setup'
   action.

\subsubsection{Julia problem files: how to use them}\label{how-julia}

   Finally, the typical use of the Julia output file is as
   follows. The {\tt PROBNAME.jl} output file is produced from the 
   {\tt PROBNAME.SIF} file by the call
   \thecode{PROBNAME = s2mpj( 'PROBNAME', injl )}
   where {\tt injl} is a struct describing the S2MPJ options (see
   Section~\ref{dec-options}) with {\tt injl.language = 'julia'} (other
   options may be specified if desired.) After including the {\tt
   PROBNAME} functions in the Julia program using 
   \thecode{include( "PROBNAME.jl" )}
   the setup of the problem is performed by a call of the form
   \thecode{pb, pbm = PROBNAME( "setup", args[:] )}
   where {\tt args[:]} is an optional comma separated list of
   problem-dependent arguments. The subsequent evaluation tasks are
   then obtained by calling one (or more) of 
   \begin{quote}\tt
       fx          = PROBNAME( "fx", pbm, x )\\
       fx, gx      = PROBNAME( "fgx", pbm, x )\\
       fx, gx, Hx  = PROBNAME( "fgHx", pbm, x )\\
       Hxv         = PROBNAME( "fHxv", pbm x, v )\\
       cx          = PROBNAME( "cx", pbm, x )\\
       cx, Jx      = PROBNAME( "cJx", pbm, x )\\
       cx, Jx, HXs = PROBNAME( "cJHx", pbm,  x )\\
       Jxv         = PROBNAME( "cJxv", pbm, x, v )\\
       cIx         = PROBNAME( "cIx", pbm, x, I )\\
       cIx, cIJx   = PROBNAME( "cIJx", pbm, x, I )\\
       cIx, cIJx, cIJHx = PROBNAME( "cIJHx", pbm,  x, I )\\
       cIHv        = PROBNAME( "cIHxv", pbm, x, v, I )\\
       Lxxy        = PROBNAME( "Lxy", pbm, x, y )\\
       Lxy, Lgxy   = PROBNAME( "Lgxy", pbm, x, y )\\
       Lxy, Lgxy, LgHxy = PROBNAME( "LgHxy", pbm, x, y )\\
       LHxyv       = PROBNAME( "LHxyv", pbm, x, y, v )\\
       LIxy        = PROBNAME( "LIxy", pbm, x, y, I )\\
       LIxy, LIgxy = PROBNAME( "LIgxy(", pbm, x, y, I )\\
       LIxy, LIgxy, LIgHxy = PROBNAME( "LIgHxy", pbm, x, y, I )\\
       LIHxyv      = PROBNAME( "LIHxyv", pbm, x, y, v, I )
   \end{quote}
   Note that {\tt pbm} always occurs as the second input argument.

\subsection{Ensuring full compatibility with the Fortran decoder}

If, for the purpose of comparison with prexisting numerical
experiments, one desires to run a version of the S2MPJ framework
which is fully compatible with the Fortran decoder, one needs to
(re)decode the SIF problems of interest with S2MPJ using an {\tt
  options} struct such that
\begin{quote}\tt
options.keepcformat = 1;\\
options.keepcorder  = 1;\\
options.pbxscale = 1;
\end{quote}
(see Section~\ref{dec-options}).

\numsection{Using reduced precision with S2MPJ in Matlab}\label{mp-s}

When using S2MPJ in the Matlab environment, it is possible to use
reduced precision aritmetic for all evaluations (of values,
derivatives and Hessian-times-vector products), provided the Matlab
Symbolic Math Toolbox is installed in the environment. If this is the
case, setting up the problem data-structures may now be done by the
call \thecode{[ pb, pbm ] = probname( 'setup\_redprec',
  problem\_parameters, ndigits ) } Note that the action keyword is now
{\tt 'setup\_redprec'} and that a new input argument {\tt
  ndigits} is now present at the end of the calling sequence.  This
new argument {\tt ndigits} specifies the number of digits to use in
the computation and it must be less or equal to 15.

If the action keyword {\tt 'setup\_redprec'} is used, the problem data
is converted to the requested number of digits at the end of the setup
process. As a consequence, its accuracy remains limited by the
double-precision computations performed during the setup operations
(and also by the accuracy of the SIF-provided reals, see the next
Section), justifiying the bound {\tt ndigits <= 15}.  Thus only the
\textit{evaluation} of function values and derivatives is affected by
the reduced-precision request. Moreover, this is only true if the
vector {\tt x} at which evaluation is requested is itself passed by
the user to the problem file in reduced-precision format. Also note
that, for now,  reduced-precision Hessians and Jacobians are no longer
computed the sparse format.

Using reduced precision is of course optional: a call to problem setup
using the action keyword {\tt 'setup'} as described in
Section~\ref{probfiles-s} results in standard double-precision
computations and sparse Hessians and Jacobians.  It is also possible
to disable the reduced-precision option when decoding a SIF problem by
setting the {\tt redprec} option of {\tt s2mpj} to 0 (see
Section~\ref{dec-options}). The resulting {\tt PROBNAME.m} file is
then marginally shorter but will issue an error message if the keyword
{\tt 'setup\_redprec'} is specified. Disabling the default of
requiring the reduced-precision option is automatic if the Symbolic
Math Toolbox is not installed.

\numsection{Testing}\label{testing-s}

The validity of the S2MPJ decoder has been tested on 1075 problems of
the \cutest\ collection (see Section~\ref{distro-s}), avoiding
large SIF files and problems involving external Fortran subroutines
(see Section~\ref{limitations-s}). The list of these problems and the
corresponding files are part of the framework (see Section~\ref{distro-s})

A first step was being able to decode all 1075 SIF files in the three
programming languages of interest. Once this was achieved, we verified
the coherence between the evaluations produced at the problem's starting
point by the S2MPJ (in standard double precision) using the Python
problem files and those produced by the Fortran decoder (as obtained
from the {\tt pycutest} wrapper). This turned out to be
subject to a significant limitation due to a feature of the Fortran
decoder inherited from the MPS format: the {\tt OUTSDIF.d} data file
written by the Fortran decoder and subsequently used by the
\cutest\ tools only stores real numbers with eight significant
digits. Unfortunately, the methods used by Fortran and Matlab (when
writing the problem files) to promote low precision numbers to
double precision ones differ: Fortran first rounds the low precision
number before completing it with zeros, while Matlab does not perform
any rounding. One therefore has to live with a coherence between the
Fortran and Python results of the order of single precision.  When
explicit numbers (such as least-squares data values) are explicitly
given in the SIF file, this can be viewed as a minor limitation. It
remains however an advantage of S2MPJ that, because the problem files
are executable, computation of constants (such as meshsizes computed
from problem dimension) not depending on SIF-provided reals are
conducted in full precision. If, as is the case for a significant
number of problems, all constants are of this type, then the relative
difference between results obtained with S2MPJ and those obtained by
the Fortran decoder typically falls below $10^{-14}$.

The final step was to verify the coherence of the evaluation produced
by the Matlab and Julia problem files with those produced by the
Python files. This was much more accurate, again typically resulting
in relative differences below $10^{-14}$.

All validation tests were run on
a Dell Precision computer with 64 GiB of memory and running Ubuntu and
Matlab R2024a, Python 3.8.10 and Julia 1.10.4.

\numsection{Writing your own optimization problem in S2PMJ format}

It is of course possible to create new test examples directly in
Matlab, Python or Julia.  For making them compatible with the S2MPJ
libraries, it is enough that they adhere to the specifications
described above (the ``S2MPJ format'').  In doing so, the choice of
detail in the GPS description of the new problem can be chosen freely:
from very detailed (typically as provided by the \cutest\ SIF files)
to very synthetic. To illustrate this point,
Appendix~\ref{rosexamples} gives several possible versions of the
ROSENBR problem (in the Matlab context and using double precision): 
\begin{itemize}
\item {\tt ROSENBR.m} is the detailed and (here too) verbose version
  produced by S2MPJ, using two
  groups of type {\tt gL2} and two nonlinear elements of type {\tt eSQ}.
\item {\tt ROSSIMP1.m} is simplified version where the detailed GPS
  structure is still used, but the (for this problem) unnecessary
  statements automatically generated by S2MPJ have been removed,
\item {\tt ROSSIMP2.m} is further simplification, using now one single
  {\tt TRIVIAL} group and two nonlinear elements of types {\tt E1} and
  {\tt E2},
\item {\tt ROSSIMP3.m} is a short synthetic version, with a single
  {\tt TRIVIAL} group and one nonlinear elements of type {\tt EROSNB}.
  This last version may be used as a template for the simplest
  possible specification of an unconstrained problem in S2PMJ format.
\end{itemize}
In addition, {\tt ROSSIMP1NE} also presents the problem as a system of
nonlinear equations (with one group per equation, the first one having
one nonlinear element of type {\tt eSQ}).

\numsection{Distribution}\label{distro-s}

The S2MPJ decoder {\tt s2mpj.m} and its associated libraries {\tt
  s2mpjlib.m}, {\tt s2mpjlib.py} and {\tt s2mpjlib.jl} are available
on the Github repository {\tt https://github.com/GrattonToint/S2MPJ}. This repository also
provides the 1075 SIF files used in our test, as well as the
corresponding Matlab, Python and Julia problems files (generated using
the S2MPJ defaults). A Matlab tool {\tt regenerate.m} is also
available to regenerate the Matlab/Python/Julia problem files (on a
Linux system), should one desire to use different S2MPJ options.

\numsection{Conclusion, discussion and perspectives}\label{concl-s}

We have presented S2PMJ a new decoder for the SIF files defining the
\cutest\ problems.  This decoder is written in Matlab and produces
executable problem files in Matlab, Python or Julia, making the
stand-alone use of the \cutest\ collection of test problems possible
directly within the Matlab, Python or Julia environments. The code for
the decoder, its associated libraries and the generated problem files
are publicly available on Github.

As it turned out, the design of S2MPJ was useful as an independent
check of the collection SIF files (a few were either
corrected or completed) and also resulted in an improved Fortran
decoder.  But the authors are well aware of the 
limitations of the S2MPJ framework:  not all
existing SIF files could be decoded, the format of the
Matlab/Python/Julia problem files can be very verbose (it has to cover
every possible SIF syntax) and they are intrinsically slower to
execute than their Fortran counterparts, as expected for JIT compiled
languages.

In particular, the 'setup' action and the first evaluation in Julia
file can occasionally be very memory-intensive and slow, to the point
of making the approach irrelevant for problems involving a large
amount of data or a large number of variables.  It is therefore clear
that \textit{both the Fortran and S2MPJ frameworks have complementary
  uses}.  Equally clear is the fact that improvements to S2MPJ
(such as automatic storage of precompiled problem files and automatic
code vectorization) are desirable.

{\footnotesize
  
\section*{Acknowledgements}
The authors wish to thank Nick Gould, Zaikun Zhang and Dominique Orban
for their useful comments, corrections and
suggestions. They are also indebted to the Polytechnic University of
Hong Kong for supporting a visit dedicated to the S2MPJ
project. Philippe Toint finally gratefully acknowledges the partial
support of the Institut National Polytechnique (INP) of Toulouse
(France).

The authors report there are no competing interests to declare. 


}

\appendix
\appnumsection{The Rosenbrock example: from automatic to
  simplified}\label{rosexamples}

\fbox{Version 0: the S2MPJ-produced file
    (two $\ell_2$ groups, two nonlinear elements)}
\begin{lstlisting}
function varargout = ROSENBR(action,varargin)
%%%%%%%%%%%%%%%%%%%%%%%%%%%%%%%%%%%%%%%%%%%%%%%%%%%%%%%%%%%%%%%%%%%%%%%%%%%
%    Problem : ROSENBR
%    The ever famous 2 variables Rosenbrock "banana valley" problem
%    Source:  problem 1 in
%    J.J. More', B.S. Garbow and K.E. Hillstrom,
%    "Testing Unconstrained Optimization Software",
%    ACM Transactions on Mathematical Software, vol. 7(1), pp. 17-41, 1981.
%    SIF input: Ph. Toint, Dec 1989.
%    classification = 'SUR2-AN-2-0'
%%%%%%%%%%%%%%%%%%%%%%%%%%%%%%%%%%%%%%%%%%%%%%%%%%%%%%%%%%%%%%%%%%%%%%%%%%%
persistent pbm;
name = 'ROSENBR';
switch(action)
    case 'setup'
        pb.name      = name;
        pbm.name     = name;
        %%%%%%%%%%%%%%%%%%%%  PREAMBLE %%%%%%%%%%%%%%%%%%%%
        v_  = configureDictionary('string','double');
        ix_ = configureDictionary('string','double');
        ig_ = configureDictionary('string','double');
        %%%%%%%%%%%%%%%%%%%%  VARIABLES %%%%%%%%%%%%%%%%%%%%
        pb.xnames = {};
        [iv,ix_] = s2mpjlib('ii','X1',ix_);
        pb.xnames{iv} = 'X1';
        [iv,ix_] = s2mpjlib('ii','X2',ix_);
        pb.xnames{iv} = 'X2';
        %%%%%%%%%%%%%%%%%%%  DATA GROUPS %%%%%%%%%%%%%%%%%%%
        pbm.A = sparse(0,0);
        [ig,ig_] = s2mpjlib('ii','G1',ig_);
        gtype{ig} = '<>';
        iv = ix_('X2');
        if(size(pbm.A,1)>=ig&&size(pbm.A,2)>=iv)
            pbm.A(ig,iv) = 1.0+pbm.A(ig,iv);
        else
            pbm.A(ig,iv) = 1.0;
        end
        pbm.gscale(ig,1) = 0.01;
        [ig,ig_] = s2mpjlib('ii','G2',ig_);
        gtype{ig} = '<>';
        iv = ix_('X1');
        if(size(pbm.A,1)>=ig&&size(pbm.A,2)>=iv)
            pbm.A(ig,iv) = 1.0+pbm.A(ig,iv);
        else
            pbm.A(ig,iv) = 1.0;
        end
        %%%%%%%%%%%%%%% GLOBAL DIMENSIONS %%%%%%%%%%%%%%%%%
        pb.n   = numEntries(ix_);
        ngrp   = numEntries(ig_);
        pbm.objgrps = [1:ngrp];
        pb.m        = 0;
        %%%%%%%%%%%%%%%%%%% CONSTANTS %%%%%%%%%%%%%%%%%%%%%
        pbm.gconst = zeros(ngrp,1);
        pbm.gconst(ig_('G2')) = 1.0;
        %%%%%%%%%%%%%%%%%%%%  BOUNDS %%%%%%%%%%%%%%%%%%%%%
        pb.xlower = -Inf*ones(pb.n,1);
        pb.xupper = +Inf*ones(pb.n,1);
        %%%%%%%%%%%%%%%%%%%% START POINT %%%%%%%%%%%%%%%%%%
        pb.x0(1:pb.n,1) = zeros(pb.n,1);
        pb.x0(ix_('X1'),1) = -1.2;
        pb.x0(ix_('X2'),1) = 1.0;
        %%%%%%%%%%%%%%%%%%%%% ELFTYPE %%%%%%%%%%%%%%%%%%%%%
        iet_ = configureDictionary('string','double');
        [it,iet_] = s2mpjlib( 'ii', 'eSQ',iet_);
        elftv{it}{1} = 'V1';
        %%%%%%%%%%%%%%%%%%% ELEMENT USES %%%%%%%%%%%%%%%%%%
        ie_ = configureDictionary('string','double');
        pbm.elftype = {};
        ielftype    = [];
        pbm.elvar   = {};
        ename = 'E1';
        [ie,ie_] = s2mpjlib('ii',ename,ie_);
        pbm.elftype{ie} = 'eSQ';
        ielftype(ie) = iet_('eSQ');
        vname = 'X1';
        [iv,ix_,pb] = s2mpjlib('nlx',vname,ix_,pb,1,[],[],[]);
        posev = find(strcmp('V1',elftv{ielftype(ie)}));
        pbm.elvar{ie}(posev) = iv;
        %%%%%%%%%%%%%%%%%%%%%% GRFTYPE %%%%%%%%%%%%%%%%%%%%
        igt_ = configureDictionary('string','double');
        [it,igt_] = s2mpjlib('ii','gL2',igt_);
        %%%%%%%%%%%%%%%%%%%% GROUP USES %%%%%%%%%%%%%%%%%%%
        [pbm.grelt{1:ngrp}] = deal(repmat([],1,ngrp));
        nlc = [];
        for ig = 1:ngrp
            pbm.grftype{ig} = 'gL2';
        end
        ig = ig_('G1');
        posel = length(pbm.grelt{ig})+1;
        pbm.grelt{ig}(posel) = ie_('E1');
        pbm.grelw{ig}(posel) = -1.0;
        %%%%%%%%%%%%%%%%%%% OBJECT BOUNDS %%%%%%%%%%%%%%%%%
        pb.objlower = 0.0;
%    Solution
% LO SOLTN                0.0
        %%%%%%%%% DEFAULT FOR MISSING SECTION(S) %%%%%%%%%%
        %%%%%% RETURN VALUES FROM THE SETUP ACTION %%%%%%%%
        pb.pbclass = 'SUR2-AN-2-0';
        varargout{1} = pb;
        varargout{2} = pbm;
% **********************
%  SET UP THE FUNCTION *
%  AND RANGE ROUTINES  *
% **********************
    %%%%%%%%%%%%%%%% NONLINEAR ELEMENTS %%%%%%%%%%%%%%%
    case 'eSQ'
        EV_  = varargin{1};
        iel_ = varargin{2};
        varargout{1} = EV_(1)*EV_(1);
        if(nargout>1)
            g_(1,1) = EV_(1)+EV_(1);
            varargout{2} = g_;
            if(nargout>2)
                H_(1,1) = 2.0;
                varargout{3} = H_;
            end
        end
    %%%%%%%%%%%%%%%%%% NONLINEAR GROUPS  %%%%%%%%%%%%%%%
    case 'gL2'
        GVAR_ = varargin{1};
        igr_  = varargin{2};
        varargout{1} = GVAR_*GVAR_;
        if(nargout>1)
            g_ = GVAR_+GVAR_;
            varargout{2} = g_;
            if(nargout>2)
                H_ = 2.0;
                varargout{3} = H_;
            end
        end
    %%%%%%%%%%%%%%%% THE MAIN ACTIONS %%%%%%%%%%%%%%%
    case {'fx','fgx','fgHx','cx','cJx','cJHx','cIx','cIJx','cIJHx','cIJxv','fHxv',...
          'cJxv','Lxy','Lgxy','LgHxy','LIxy','LIgxy','LIgHxy','LHxyv','LIHxyv'}
        if(isfield(pbm,'name')&&strcmp(pbm.name,name))
            pbm.has_globs = [0,0];
            [varargout{1:max(1,nargout)}] = s2mpjlib(action,pbm,varargin{:});
        else
            disp(['ERROR: please run ',name,' with action = setup'])
            [varargout{1:nargout}] = deal(repmat(NaN,1:nargout));
        end
    otherwise
        disp([' ERROR: unknown action ',action,' requested from ',name,'.m'])
    end
return
end
%%%%%%%%%%%%%%%%%%%%%%%%%%%%%%%%%%%%%%%%%%%%%%%%%%%%
\end{lstlisting}

\vspace*{2mm}
\noindent
\fbox{Version 1: the uncluttered full GPS description
  (two $\ell_2$ groups, two nonlinear elements)}
\begin{lstlisting}
function varargout = ROSSIMP1(action,varargin)
%%%%%%%%%%%%%%%%%%%%%%%%%%%%%%%%%%%%%%%%%%%%%%%%%%%%%%%%%%%%%%%%%%%%%%%%%%%
%    Problem : ROSSIMP1
%    The ever famous 2 variables Rosenbrock "banana valley" problem.
%    This version is constructed from ROSENBR.m, removing all (here)
%    unnecessary statements automatically generated by S2MPJ.
%    classification = 'SUR2-AN-2-0'
%%%%%%%%%%%%%%%%%%%%%%%%%%%%%%%%%%%%%%%%%%%%%%%%%%%%%%%%%%%%%%%%%%%%%%%%%%%
persistent pbm;
name = 'ROSSIMP1';
switch(action)
    case 'setup'
        pb.name      = name;
        pbm.name     = name;
        %%%%%%%%%%%%%%%%%%%  DATA GROUPS %%%%%%%%%%%%%%%%%%%
        pbm.A = sparse(2,2);
        pbm.A(1,2) = 1.0;
        pbm.gscale(1) = 0.01;
        pbm.A(2,1) = 1.0;
        %%%%%%%%%%%%%%% GLOBAL DIMENSIONS %%%%%%%%%%%%%%%%%
        pb.n   = 2;
        ngrp   = 2;
        pbm.objgrps = [ 1, 2 ];
        pb.m        = 0;
        %%%%%%%%%%%%%%%%%%% CONSTANTS %%%%%%%%%%%%%%%%%%%%%
        pbm.gconst =[ 0.0; 1 ];
        %%%%%%%%%%%%%%%%%%%%  BOUNDS %%%%%%%%%%%%%%%%%%%%%
        pb.xlower = [ -Inf; -Inf ];
        pb.xupper = [ +Inf; +Inf ];
        %%%%%%%%%%%%%%%%%%%% START POINT %%%%%%%%%%%%%%%%%%
        pb.x0 = [ -1.2; 1.0 ];
        %%%%%%%%%%%%%%%%%%% ELEMENT USES %%%%%%%%%%%%%%%%%%
        pbm.elftype{1} = 'eSQ';
        pbm.elvar{1}   = [ 1 ];
        %%%%%%%%%%%%%%%%%%%% GROUP USES %%%%%%%%%%%%%%%%%%%
        pbm.grftype = { 'gL2', 'gL2' };
        pbm.grelt{1} = [ 1 ];
        pbm.grelw{1} = [ -1.0 ];
        %%%%%%%%%%%%%%%%%%% OBJECT BOUNDS %%%%%%%%%%%%%%%%%
        pb.objlower = 0.0;
        %%%%%%%%% DEFAULT FOR MISSING SECTION(S) %%%%%%%%%%
        %%%%%% RETURN VALUES FROM THE SETUP ACTION %%%%%%%%
        pb.pbclass = 'SUR2-AN-2-0';
        varargout{1} = pb;
        varargout{2} = pbm;
    %%%%%%%%%%%%%%%% NONLINEAR ELEMENTS %%%%%%%%%%%%%%%
    case 'eSQ'
        EV_  = varargin{1};
        varargout{1} = EV_(1)*EV_(1);
        if(nargout>1)
            g_(1,1) = EV_(1)+EV_(1);
            varargout{2} = g_;
            if(nargout>2)
                H_(1,1) = 2.0;
                varargout{3} = H_;
            end
        end
    %%%%%%%%%%%%%%%%%% NONLINEAR GROUPS  %%%%%%%%%%%%%%%
    case 'gL2'
        GVAR_ = varargin{1};
        varargout{1} = GVAR_*GVAR_;
        if(nargout>1)
            g_ = GVAR_+GVAR_;
            varargout{2} = g_;
            if(nargout>2)
                H_ = 2.0;
                varargout{3} = H_;
            end
        end
    %%%%%%%%%%%%%%%% THE MAIN ACTIONS %%%%%%%%%%%%%%%
    case {'fx','fgx','fgHx','cx','cJx','cJHx','cIx','cIJx','cIJHx','cIJxv','fHxv',...
          'cJxv','Lxy','Lgxy','LgHxy','LIxy','LIgxy','LIgHxy','LHxyv','LIHxyv'}
        if(isfield(pbm,'name')&&strcmp(pbm.name,name))
            pbm.has_globs = [0,0];
            [varargout{1:max(1,nargout)}] = s2mpjlib(action,pbm,varargin{:});
        else
            disp(['ERROR: please run ',name,' with action = setup'])
            [varargout{1:nargout}] = deal(repmat(NaN,1:nargout));
        end
    otherwise
        disp([' ERROR: unknown action ',action,' requested from ',name,'.m'])
    end
return
end
%%%%%%%%%%%%%%%%%%%%%%%%%%%%%%%%%%%%%%%%%%%%%%%%%%%%


\end{lstlisting}

\vspace*{2mm}
\noindent
\fbox{Version 2:  an intermediate version ((one trivial group, two
  nonlinear elements)}
\begin{lstlisting}
function varargout = ROSSIMP2(action,varargin)
%%%%%%%%%%%%%%%%%%%%%%%%%%%%%%%%%%%%%%%%%%%%%%%%%%%%%%%%%%%%%%%%%%%%%%%%%%%
%    Problem : ROSSIMP2
%    The ever famous 2 variables Rosenbrock "banana valley" problem
%    This version is constructed from ROSSIMP1.m, reducing from 2 nonlinear
%    groups and 1 nonlinear element to 1 trivial group and 2 nonlinear elements.
%    S2MPJ input: Ph. Toint, June 2024
%    classification = 'SUR2-AN-2-0'
%%%%%%%%%%%%%%%%%%%%%%%%%%%%%%%%%%%%%%%%%%%%%%%%%%%%%%%%%%%%%%%%%%%%%%%%%%%
persistent pbm;
name = 'ROSSIMP2';
switch(action)
    case 'setup'
        pb.name     = name;
        pbm.name    = name;
        pb.n = 2;
        pbm.objgrps = [ 1 ];
        pb.m = 0;
        pb.xlower = [ -Inf; -Inf ];
        pb.xupper = [ +Inf; +Inf ];
        pb.x0 = [ -1.2,; 1.0 ];
        pbm.elftype = { 'E1', 'E2'};
        pbm.elvar   = { [ 1, 2 ],  [ 1 ] };
        pbm.grelt{1} = [ 1, 2 ];
        pb.objlower = 0.0;
        pb.pbclass = 'SUR2-AN-2-0';
        varargout{1} = pb;
        varargout{2} = pbm;
    case 'E1'
        x  = varargin{1};
        varargout{1} = 100.0*(x(2)-x(1)*x(1))^2;
        if(nargout>1)
            g_(1,1) = -400.0*(x(2)-x(1)*x(1))*x(1);
            g_(2,1) = 200.0*(x(2)-x(1)*x(1));
            varargout{2} = g_;
            if(nargout>2)
                H_(1,1) = 1200.0*x(1)*x(1);
                H_(1,2) = -400.0;
                H_(2,2) =  200.0;
                varargout{3} = H_;
            end
        end
    case 'E2'
        x  = varargin{1};
        varargout{1} = (x(1)-1.0)^2;
        if(nargout>1)
            g_(1,1) = 2*(x(1)-1.0);
            varargout{2} = g_;
            if(nargout>2)
                H_(1,1) = 2.0;
                varargout{3} = H_;
            end
        end
    %%%%%%%%%%%%%%%% THE MAIN ACTIONS %%%%%%%%%%%%%%%
    case {'fx','fgx','fgHx','cx','cJx','cJHx','cIx','cIJx','cIJHx','cIJxv','fHxv',...
          'cJxv','Lxy','Lgxy','LgHxy','LIxy','LIgxy','LIgHxy','LHxyv','LIHxyv'}
        if(isfield(pbm,'name')&&strcmp(pbm.name,name))
            pbm.has_globs = [0,0];
            [varargout{1:max(1,nargout)}] = s2mpjlib(action,pbm,varargin{:});
        else
            disp(['ERROR: please run ',name,' with action = setup'])
            [varargout{1:nargout}] = deal(repmat(NaN,1:nargout));
        end
    otherwise
        disp([' ERROR: unknown action ',action,' requested from ',name,'.m'])
    end
return
end
%%%%%%%%%%%%%%%%%%%%%%%%%%%%%%%%%%%%%%%%%%%%%%%%%%%%


\end{lstlisting}

\vspace*{2mm}
\noindent
\fbox{Version 3: the short synthetic version (one trivial group, one
  nonlinear element)}
\begin{lstlisting}
function varargout = ROSSIMP3(action,varargin)
%%%%%%%%%%%%%%%%%%%%%%%%%%%%%%%%%%%%%%%%%%%%%%%%%%%%%%%%%%%%%%%%%%%%%%%%%%%
%    Problem : ROSSIMP3
%    The ever famous 2 variables Rosenbrock "banana valley" problem
%    This version uses  1 trivial group and 2 nonlinear elements.
%    classification = 'SUR2-AN-2-0'
%%%%%%%%%%%%%%%%%%%%%%%%%%%%%%%%%%%%%%%%%%%%%%%%%%%%%%%%%%%%%%%%%%%%%%%%%%%
persistent pbm;
name = 'ROSSIMP3';
switch(action)
    case 'setup'
        pb.name      = name;
        pbm.name     = name;
        pb.n = 2;
        pb.m = 0;
        pbm.objgrps  = [ 1 ];
        pb.xlower    = [ -Inf; -Inf ];
        pb.xupper    = [ +Inf; +Inf ];
        pb.x0        = [ -1.2,; 1.0 ];
        pbm.elftype  = { 'EROSNB'};
        pbm.elvar    = { [ 1, 2 ] };
        pbm.grelt{1} = [ 1 ];
        pb.objlower  = 0.0;
        pb.pbclass   = 'SUR2-AN-2-0';
        varargout{1} = pb;
        varargout{2} = pbm;
    case 'EROSNB'    %  The ever famous "banana" function!
        x  = varargin{1};
        varargout{1} = 100.0*(x(2)-x(1)*x(1))^2+(x(1)-1.0)^2;
        if(nargout>1)
            g_(1,1) = -400.0*(x(2)-x(1)*x(1))*x(1)+2.0*(x(1)-1.0);
            g_(2,1) = 200.0*(x(2)-x(1)*x(1));
            varargout{2} = g_;
            if(nargout>2)
                H_(1,1) = 1200.0*x(1)*x(1)+2.0;
                H_(1,2) = -400.0;
                H_(2,2) =  200.0;
                varargout{3} = H_;
            end
        end
    %%%%%%%%%%%%%%%% THE MAIN ACTIONS %%%%%%%%%%%%%%%
    case {'fx','fgx','fgHx','cx','cJx','cJHx','cIx','cIJx','cIJHx','cIJxv','fHxv',...
          'cJxv','Lxy','Lgxy','LgHxy','LIxy','LIgxy','LIgHxy','LHxyv','LIHxyv'}
        if(isfield(pbm,'name')&&strcmp(pbm.name,name))
            pbm.has_globs = [0,0];
            [varargout{1:max(1,nargout)}] = s2mpjlib(action,pbm,varargin{:});
        else
            disp(['ERROR: please run ',name,' with action = setup'])
            [varargout{1:nargout}] = deal(repmat(NaN,1:nargout));
        end
    otherwise
        disp([' ERROR: unknown action ',action,' requested from ',name,'.m'])
    end
return
end
%%%%%%%%%%%%%%%%%%%%%%%%%%%%%%%%%%%%%%%%%%%%%%%%%%%%


\end{lstlisting}

\vspace*{2mm}
\noindent
\fbox{Version 4: as a sytem of two nonlinear equations}
\begin{lstlisting}
function varargout = ROSSIMP1NE(action,varargin)
%%%%%%%%%%%%%%%%%%%%%%%%%%%%%%%%%%%%%%%%%%%%%%%%%%%%%%%%%%%%%%%%%%%%%%%%%%%
%    Problem : ROSSIMP1NE
%    The ever famous 2 variables Rosenbrock "banana valley" problem
%    This version considers the problem as a system of nonlinear equations.
%    classification = 'NOR2-AN-2-0'
%%%%%%%%%%%%%%%%%%%%%%%%%%%%%%%%%%%%%%%%%%%%%%%%%%%%%%%%%%%%%%%%%%%%%%%%%%%
persistent pbm;
name = 'ROSSIMP1NE';
switch(action)
    case 'setup'
        pb.name      = name;
        pbm.name     = name;
        %%%%%%%%%%%%%%%%%%%  DATA GROUPS %%%%%%%%%%%%%%%%%%%
        pbm.A = sparse(2,2);
        pbm.A(1,2) = 1.0;
        pbm.gscale(1) = 0.1;
        pbm.A(2,1) = 1.0;
        %%%%%%%%%%%%%%% GLOBAL DIMENSIONS %%%%%%%%%%%%%%%%%
        pb.n   = 2;
        pbm.congrps = [ 1, 2 ];
        pb.neq      = 2;
        pb.m        = 2;
        %%%%%%%%%%%%%%%%%%% CONSTANTS %%%%%%%%%%%%%%%%%%%%%
        pbm.gconst =[ 0.0; 1 ];
        %%%%%%%%%%%%%%%%%%%%  BOUNDS %%%%%%%%%%%%%%%%%%%%%
        pb.xlower = [ -Inf; -Inf ];
        pb.xupper = [ +Inf; +Inf ];
        %%%%%%%%%%%%%%%%%%%% START POINT %%%%%%%%%%%%%%%%%%
        pb.x0 = [ -1.2; 1.0 ];
        %%%%%%%%%%%%%%%%%%% ELEMENT USES %%%%%%%%%%%%%%%%%%
        pbm.elftype{1} = 'eSQ';
        pbm.elvar{1}   = [ 1 ];
        %%%%%%%%%%%%%%%%%%%% GROUP USES %%%%%%%%%%%%%%%%%%%
        pbm.grelt{1} = [ 1 ];
        pbm.grelw{1} = [ -1.0 ];
        %%%%%%%%%%%%%%%%%%% OBJECT BOUNDS %%%%%%%%%%%%%%%%%
        pb.objlower = 0.0;
        %%%%%% RETURN VALUES FROM THE SETUP ACTION %%%%%%%%
        pb.pbclass = 'NOR2-AN-2-0';
        varargout{1} = pb;
        varargout{2} = pbm;
    case 'eSQ'
        EV_  = varargin{1};
        varargout{1} = EV_(1)*EV_(1);
        if(nargout>1)
            g_(1,1) = EV_(1)+EV_(1);
            varargout{2} = g_;
            if(nargout>2)
                H_(1,1) = 2.0;
                varargout{3} = H_;
            end
        end
    %%%%%%%%%%%%%%%% THE MAIN ACTIONS %%%%%%%%%%%%%%%
    case {'fx','fgx','fgHx','cx','cJx','cJHx','cIx','cIJx','cIJHx','cIJxv','fHxv',...
          'cJxv','Lxy','Lgxy','LgHxy','LIxy','LIgxy','LIgHxy','LHxyv','LIHxyv'}
        if(isfield(pbm,'name')&&strcmp(pbm.name,name))
            pbm.has_globs = [0,0];
            [varargout{1:max(1,nargout)}] = s2mpjlib(action,pbm,varargin{:});
        else
            disp(['ERROR: please run ',name,' with action = setup'])
            [varargout{1:nargout}] = deal(repmat(NaN,1:nargout));
        end
    otherwise
        disp([' ERROR: unknown action ',action,' requested from ',name,'.m'])
    end
return
end
%%%%%%%%%%%%%%%%%%%%%%%%%%%%%%%%%%%%%%%%%%%%%%%%%%%%

\end{lstlisting}

\end{document}